\documentclass{article}

\begin{document}

Tsemo Aristide

Visitor, University of Toronto

Department of mathematics

100 St Georges Street

tsemoaristide@hotmail.com

\bigskip

\bigskip

             {\bf THE DIFFERENTIAL GEOMETRY OF COMPOSITION
             SEQUENCES OF DIFFERENTIABLE MANIFOLDS.}

\bigskip

\bigskip

\centerline{\bf Abstract.}

{\it Let $F_0=B,F_1,..,F_n$ be a sequence of differentiable
manifolds and $G_l, l\geq 1$ a Lie subgroup of the group of
diffeomorphisms of $F_l$, $H_l$ a central subgroup of $G_l$. We
denote by $K_l$ the quotient of $G_l$ by $H_l$.  We suppose also
given a locally trivial principal bundle $p_{K_l}$ over $F_{l-1}$
which typical fiber is $K_l$. In this paper we study the
differential geometry of the problem of the existence of a
composition sequence of manifolds $M_n\rightarrow
M_{n-1}\rightarrow..\rightarrow M_0=B$, where each map
$f_l:M_{l}\rightarrow M_{l-1}$ is a locally trivial differentiable
bundle over $M_{l-1}$ which typical fiber is $F_l$, and such that
the image of the transitions functions is contained in $G_l$. We
associated to this problem a tower of gerbes, and we study its
differential geometry: that is we define the notion of connective
structure and curvature.}

\bigskip

\centerline{\bf Introduction.}

\bigskip

Let $F_0=B,F_1,..,F_n$ be a sequence of differentiable manifolds
and $G_l, l\geq 1$ a Lie subgroup of the group of diffeomorphisms
of $F_l$, $H_l$ a central subgroup of $G_l$. We denote by $K_l$
the quotient of $G_l$ by $H_l$.  We suppose also given a principal
bundle $p_{K_l}$ over $F_{l-1}$ which typical fiber is $K_l$. In
this paper we study the differential geometry of the problem of
the existence of a composition sequence of manifolds
$M_n\rightarrow M_{n-1}\rightarrow..\rightarrow M_0=B$, where each
map $f_l:M_{l}\rightarrow M_{l-1}$ is a locally trivial
differentiable bundle over $M_{l-1}$ which typical fiber is $F_l$,
and such that the image of the transitions functions is contained
in $G_k$. We can associate to $f_l$ a bundle over $M_{l-1}$ which
typical fiber is $G_l$, we suppose also that the quotient of the
restriction of this bundle to a fiber of $f_{l-1}$ by $H_l$ is
$p_{K_l}$. We associated to this problem a tower of gerbes, and we
study its differential geometry.

The construction of those sequences of fibrations fits in the
general problem of the construction of structures in mathematics.
There are essentially two ways to obtain new structures:

One is the completion procedure, which is used generally to create
space in which unsolved problems have a solution. It  also shows
the link between  algebra and analysis.

The other essential manner to construct structures in mathematics
is the gluing procedure, which is the bridge from analysis to
geometry. It is the one involved in our problem. To solve the
gluing problem, that is to find obstructions to define global
objects which are defined locally, one needs cohomology theories
interpreted geometrically.

More practically, consider a topological manifold $M$, it is built
by gluing open subsets of a vector space say $E$, to study $M$, we
can for example generalize tools used to study $E$. For instance,
we can define the notion of continuous functions. Can functions
defined locally  be defined globally ? This problem is solved
using a $1-$ cohomology theory.

The notion of sheaf of functions leads to the notion of locally
trivial bundles. which are  particular case of sheaves of
functions. One can also try to solve the following problem. Can
sheaves of bundles defined locally  be defined globally? To solve
this problem, one needs the notion of stacks defined by Giraud
[4]. The obstructions to solve this problem is given by a
$2-$cohomology theory.

More generally, one can try to define the notion of sheaf of
stacks ($3-$sheaf) which  would lead to the notion of a
$3-$cohomology theories,..., the notion of $n-$sheaf, which would
lead to the notion of a $n-$dimension cohomology theory.

 To define  a notion of $n-$sheaf, one needs to define
 first, the notion of $n-$category, which is not well understood
 nowadays.

The gluing procedure  in the $1-$dimensional case  is the theory
of sheaves. It is used to construct Hilbert schemes (see Gottsche
[5]).

In the two dimensional case, it has been established by Giraud ,
it is the theory of stacks.

\bigskip

{\bf I. The first lifting problem.}

\medskip

 Let $B$ and $F$ be two differentiable manifolds, $ G$ a Lie
subgroup of diffeomorphisms of $F$, $H$ a subgroup of $G$ central
in $G$, we denote by $K$ the quotient $G/H$. We suppose defined on
$B$ a locally trivial $K-$principal bundle. We want to study the
differential geometry of the following problem: Classify the
locally trivial bundles over $B$, which typical fiber is $F$ and
such that the family of elements $g_{ij}$ which define each
trivialization $(U_i,g_{ij},i,j \in I)$  is included in $G$. This
means that the $Diff(F)$ bundle associated to the bundle has a
$G-$reduction. We can associate also to a solution of our problem
a principal bundle $p_G$ over $B$ which typical fiber is $G$. We
suppose moreover that the quotient of $p_G$ by $H$ is $p_K$. This
problem will be called the first lifting problem.

\medskip

 We have the exact sequence $1\rightarrow H\rightarrow
G\rightarrow K \rightarrow 1$. We will assume that the map
$G\rightarrow K$ has local sections.

\medskip

Our problem is equivalent to the following: study the existence of
a bundle over $B$ which typical fiber is $G$, and which quotient
by $H$ is  $p_K$. To solve it, we will define a gerbe over $B$
which classifying cocycle is the obstruction of the existence of
such a bundle. This problem will be called the principal problem
associated to the first lifting problem.

Let's recall now some facts from gerbe theory.

\medskip

{\bf Definition 1.}

Let $B$ be a manifold, a sheaf $S$ of categories on $B$, is a map
$U\rightarrow S(U)$, where $U$ is an open set of $B$, and $S(U)$
category which satisfies the following properties:

- To each inclusion $U\rightarrow V$, there exists a map
$r_{U,V}:S(V)\rightarrow S(U)$ such that $r_{U,V}\circ
r_{V,W}=r_{U,W}$.

- Gluing conditions for objects,

Consider  a covering family $(U_i)_{i\in I}$ of an open set $U$ of
$B$, and for each $i$, an object $x_i$ of $S(U_i)$, suppose that
there exists a map $g_{ij}:r_{U_i\cap U_j,U_j}(x_j)\rightarrow
r_{U_i\cap U_j, U_i}(x_i)$ such that $g_{ij}g_{jk}=g_{ik}$, then
there exists an object $x$ of $C(U)$ such that $r_{U_i,U}(x)=x_i$

Gluing conditions for arrows,

Consider two objects $P$ and $Q$ of $S(B)$, then the map
$U\rightarrow Hom(r_{U,B}(P),r_{U.B}(Q))$ is a sheaf.

Moreover, if the following conditions are satisfied the sheaf of
categories $S$ is called a gerbe

$G1$

There exists a covering family $(U_i)_{i\in I}$ of $B$ such that
for each $i$ the category $S(U_i)$ is not empty

$G2$

Let $U$ be an open set of $B$, for each objects $x$ and $y$ of
$U$, there exists a covering $(U_i)_{i\in I}$ of $U$ such that
$r_{U_i,U}(x)$ and $r_{U_i,U}(y)$ are isomorphic

$G3$

Every arrow of $S(U)$ is invertible, and there exists a sheaf $A$
in groups on $B$, such that for each object $x$ of $S(U)$,
$Hom(x,x)= A(U)$, an this family of isomorphisms commute with the
restriction
 maps.

The sheaf $A$ is called the band of the gerbe $S$, in the sequel,
we will consider only gerbes with commutative band.

\medskip

{\bf Notation.}

\medskip

 for a covering family $(U_i)_{i\in I}$ of $B$, and an object $x_i$
 of $S(U_i)$, we denote
by $x^i_{i_1..i_n}$ the element $r_{({U_{i_1}}\cap..\cap
U_{i_n},U_i)}(x_i)$, and by $U_{i_1..i_n}$ the intersection
$U_{i_1}\cap..\cap U_{i_n}$.

\medskip

Let endows $B$ with a gerbe $S$ with band $A$, one can associates
the following two Cech cocycle to $B$. Consider a covering family
$(U_i)_{i\in I}$ such that for each $i$, $S(U_i)$ is not empty,
let $x_i$ and $x_j$ be respectively objects of $S(U_i)$ and
$S(U_j)$, we consider a map $g_{ij}: {{x}^j}_{ij}\rightarrow
{{x}^i}_{ij}$, the map $t_{ijk}=g_{ij}g_{jk}g_{ki}$ of
$Aut({{{x}^i}_{ijk}})$ is the classifying $2$ Cech cocycle.

\medskip

{\bf Theorem 1. see [4]}

{\it Two gerbes defined over $B$ with band $A$ are isomorphic if
and only if their  associated  $2-$cocycles define the same
cohomological class. Conversely for each  $2-$Cech $A$ cocycle
$c$, one can define a gerbe which associated cocycle is $c$.}

\bigskip

We define now on $B$ the sheaf of categories used to solve our
problem. For each open set $U$ of $B$, we denote by $C(U)$ the
category which objects are local trivial bundles over $U$ which
typical fiber is $F$, such that the $Diff(F)$ bundle associated
has a $G-$reduction. This means that the image of the transitions
functions are elements of $G$. Moreover we suppose that the
quotient of this reduction by $H$ is the restriction ${p_K}_U$, of
$p_K$ to $U$.
 A map between two objects of $C(U)$, is an isomorphism of bundle
 which induces a map between their associated principal bundles
 which pushes forward to the identity of ${p_K}_U$. This means
 that it is an isomorphism which preserves and acts on each fiber
 by an element of $H$.

The map
$$
U\longrightarrow C(U)
$$
defines on $B$ a gerbe which band is $H'$, the sheaf of
differentiable $H$ functions defined on $B$.

  Let
  $$
  k_{ij}: U_i\cap U_j\longrightarrow K
  $$
be the the transitions functions of the bundle $p_K$ associated to
the covering family $(U_i)_{i\in I}$.

Consider a map
$$
g_{ij}:U_i\cap U_j\longrightarrow G
$$
which pushes forward to $k_{ij}$. The map $g_{ij}g_{jk}g_{ki}$
defines the classifying cocycle $c_H$ of the gerbe $C$.

\medskip

We can also defined the principal gerbe $C_P$ associated to our
problem, for each open set $U$ of $B$, we will denote by $C_P(U)$,
the set of locally trivial bundles defined on $U$ which typical
fiber is $G$, and which pushes forward to the restriction of $p_K$
to $U$.

\bigskip

{\bf Connections on principal bundles.}

\bigskip

Let ${\cal K}$, ${\cal G}$ and ${\cal H}$ be respectively the Lie
algebra of $K$, $G$ and $H$. A connection on $p_K$ is a $1-$ form
$\omega: TB\rightarrow {\cal K}$ such that:

1

 For each fundamental vector field $A_B$ defined by the element
$A$ of ${\cal K}$, we have: $\omega(A_B(x))=A$,

2

For each element $g$ of $K$, we have
$\omega(g^*(X))=Ad(g^{-1})(X))$.

\medskip

The curvature of $\omega$ is the $2-$form $d\omega+{1\over
2}[\omega,\omega]$.

\medskip

{\bf The notion connective structure of the principal gerbe.}

\medskip

Let $\omega$ be a connection on $p_K$, and $e_U$  an element of
$C_P(U)$, recall that is  a principal bundle  over $U$ which
typical fiber is $G$. We denote by $Co(e_U,\omega)$ the family of
connections defined on $e_U$ which push forward to the restriction
of $\omega$ to $U$. Consider the map $p_1:{\cal G}\rightarrow{\cal
K}$, and the map $p_2:e_U\rightarrow p_{K,U}$ where $p_{K,U}$ is
the restriction of $p_K$ to $U$. the fact that the connection $a$
defined on $e_U$ pushes forward to the restriction $\omega_U$ of
$\omega$ to $U$ means that ${p_2}^*(w_U)=p_1(a)$.

 The map $U\longrightarrow Co(e_U,\omega)$ is
called the connective structure associated to $\omega$.

\medskip

Supposed that $((U_i)_{i\in I},g_{ij})$ is a trivialization of the
bundle $p_K$. Let $e_i$ be an element of $C_P(U_i)$ endowed with a
connection $w_i$ over the restriction of $\omega$ to
$(p_K)_{U_i}$, the restriction of $p_K$ to $U_i$.

There exists a map $g_{ij}$ between $e^j_{ij}$ and $e^i_{ij}$ the
respective restrictions of $e_j$ and $e_i$ to $U_{ij}$,

Consider the form $w_j-g_{ij}^*w_i$,  (it pushes forward  to a
form $a_{ij}$ defined on $TU_{ij}$ since $H$ is central in $G$),
 and ${e^i}_{jk}$ the restriction of $e_i$ to $U_{ijk}$. On
$U_{ijk}$, we have the restriction of the forms $w_j-g_{ij}^*w_i$,
$w_k-g_{jk}^*w_j$ and $w_k-g_{ik}^*w_i$, define respectively on
${e^j}_{ik}$, ${e^k}_{ij}$, we can define the form

$$
g_{jk}^*(w_j-g_{ij}^*w_i)+(w_k-g_{jk}^*w_j)-(w_k-g_{ik}^*w_i)=
$$

$$
-(g_{ij}g_{jk})^*w_i+g_{ik}^*w_i=g_{ik}^*(w_i-(g_{ij}g_{jk}g_{ki})^*w_i)=
g_{ik}^*(w_i-c_{ijk}^*w_i)=c_{ijk}^{-1}dc_{ijk}
$$

\bigskip

{\bf Proposition 1.} {\it The family $d(a_{ij})$ is a
$1-$cocycle.}

\medskip

{\bf Proof.}

On $U_{ijk}$, we have $a_{jk}-a_{ik}+a_{ij}=c_{ijk}^{-1}dc_{ijk}$,
since $d(c_{ijk}^{-1}dc_{ijk})=0$, we obtain the result.

\bigskip

Consider now the curvature $K(w_i)$ of the connection $w_i$, it is
the form $w_i+{1\over 2}[w_i,w_i]$. For any other form $w'_i$ over
$\omega_i$, let $a_i=w_i-w'_i$, we have $K(w_i)=K(w'_i)+da_i$,
this implies that $dK(w_i)$ is does not depend of the form over
$\omega_i$.
 To set in a conceptual theory, one as to define
the cotangent and tangent spaces of the principal gerbe.

\medskip

{\bf Definition 2.}

Let $U$ be an open set of $B$,  an element of the tangent space of
$C(U)$ will be a family of vectors fields $(X_e)$ where $e$ is an
object of $C(U)$, such that if $e$ and $f$ are objects of $C(U)$,
then there exists a map $g:e\rightarrow f$ such that
$Tg(X_e)=X_f$. An element of the cotangent of will be a family of
forms $(\alpha)_e$ such that for each vector field $(X_e)$, we
have $\alpha_f(X_f)=\alpha_e(X_e)$, one can define in an obvious
way $\Lambda^*C$.

A $n-$tangent field of $C(U)$ $nTC(U)$ will be a family of of
$n-$vector fields $((X^1_e,...,X^n_e) e\in Obj(C(U))$ where
$X^i_e$ is a vector field tangent to the element $e$ of $C(U)$
such for each object $e$ and $f$ of $C(U)$, there is a map
$g:e\rightarrow f$ such that $Tg(X^i_e)=X^i_f$.

\medskip

 one deduce that the family of
$3-$form $dK(w_i)$ define a form $\Omega$ on $3TC$ which takes
value in ${\cal G}$ this form is called the curvature of the
associated gerbe.

\bigskip

{\bf II. The generalization.}

\medskip

Let $B=F_{0}, F_1,...,F_n$ be a sequence of manifolds. For each
$F_i$, we consider a Lie subgroup $G_i$ of diffeomorphisms  of
$F_i$, and a central subgroup $H_i$ of $G_i$. Let denote by $K_i$
the quotient of $G_i$ by $H_i$, we also suppose given locally
trivial principal bundles $p_{K_i}$ over $F_{i-1}$ which typical
fiber is $K_i$. We will denote by ${\cal H}_i$, ${\cal G}_i$ and
${\cal K}_i$ the respective Lie algebras of $H_i$, $G_i$ and
$K_i$.

We want to study the differentiable geometry of the following
problem:

classify sequences $M_n\rightarrow
M_{n-1}\rightarrow..M_1\rightarrow M_0=B$, where each map
$f_i:M_{i+1}\rightarrow M_i$ is locally trivial differentiable
bundle which typical fiber is $F_{i+1}$, such that the transitions
functions
 are elements of $G_{i+1}$,  so we can define a locally trivial principal
 bundle $p_{{MGii+1}}$ over $M_i$ which typical fiber is $G_{i+1}$,
 and consider is quotient $p_{MKii+1}$ by $H_{i+1}$.
 We suppose that the restriction $p_{FKii+1}$ of $p_{MGii+1}$ to
 a fiber of $f_{i-1}$ (naturally diffeomorphic to $F_i$) is
 $p_{K_{i+1}}$. We will call this problem the $n-$lifting problem.

 \bigskip

 We will associate to our problem, a tower of gerbes, and cocycles
 which will be the obstructions to solve it.

\medskip

First we define on $B$ the gerbe $C_1$ associated to the first
lifting problem defined by the data $F_1$, $H_1$, $G_1$ and $K_1$.

Let $U$ be an open set of $B$ and $e_1(U)$ an element of $C_1(U)$,
we can associate to $e_1(U)$, the gerbe $C_2(e_1)$ which solve the
first lifting problem which data are $F_2$, $G_2$ and $K_2$, over
the basis $e_1(U)$. Remark that $p_{K_2}$ induces naturally a
$K_2-$principal bundle $p_{e_1}$ over $e_1(U)$. The quotient of
solutions the principal $G_2-$bundles associated to the solutions
of this first lifting problem by $H_2$ must be restriction of
$p_{e_1}$.

Suppose that we have defined the gerbes $C_1$,
$C_2(e_1)$,..,$C_i(e_{i-1})$, let $e_i$ be an object of
$C_i(e_{i-1})$, we define the gerbe $C_{i+1}(e_i)$ which solve the
first lifting problem associated to $F_{i+1}$, $G_{i+1}$,
$K_{i+1}$ with basis $e_i$. Remark that $p_{K_{i+1}}$ induces
naturally a $K_{i+1}-$principal bundle over $e_i$, we denote it by
$p_{e_i}$ . We suppose  that  quotient of the principal $G_{i+1}$
bundles associated to the solutions of this first lifting problem
by $H_{i+1}$ are restriction of $p_{e_i}$.

\medskip

 We have assumed that the group $H_i$ is central in $G_i$,
 let $(U^i_{jk},{k^i}_{jk})$ be the trivialization of the bundle $p_{K_i}$,
 one defines over $F_{i-1}$ a locally trivial bundle $Lie_{ii-1}$ which typical
 fiber is $H_i$ and which trivialization maps are:
 $$
 U_{jk}\times {\cal G}_i\longrightarrow U_{jk}\times {\cal G}_i
 $$

$$
(x,y)\longrightarrow (x,Ad({g^i}_{jk})(y))
$$
where ${g^i}_{jk}$ is an element of $G_i$ over ${k^i}_{jk}$.

\medskip

this sheaf $Lie_{ii-1}$ is defined for the first lifting problem
associated to $F_i$, $H_i$, $G_i$ and $K_i$ which basis space is
$F_{i-1}$.

Suppose that we have defined a sheaf over $B$,
$Lie_{ii-1..1}(F_i,H_i,G_i)$ associated to the $i-$lifting problem
defined by the family of groups $F_1$,..,$F_{i}$,  $H_1,..,H_i$,
$G_1,..,G_i$. We can define the $i-$lifting problem associated to
the family $F'_1=F_1$,..,$F'_{i-1}=F_{i-1}$, $F'_i= Lie_{ii-1}$
$H'_1=H_1,..,H'_{i-1}=H_{i-1}, H'_i$, $G_1,..,G_{i-1}, G'_i$,
where $G'_i$ and $H'_i$ are  the groups of automorphisms of
$Lie_{ii-1}$ which elements project respectively onto elements of
$G_i$ and $H_i$. We set $Lie_{i+1..1}(F_{i+1},G_{i+1},H_{i+1})=
Lie_{i..1}({F'_i},G'_i,H'_i)$.

\bigskip

{\bf The classifying cocycle.}

\medskip

We will now define a $n+1-$cocycle associated to the tower of
gerbes.

\medskip

For each manifold $F_i$, we will consider a trivialization
$(({U^i}_j,{ k^i}_{j_1j_2}, j\in J)$ of the bundle $p_{K_i}$.
Recall that $B=F_0$.

For each open set ${U^0}_{j_1}$, consider an object ${e^1}_{j_1}$
of $C_1({U^0}_{j_1})$. There exists a map ${g^1}_{j_1j_2}$ over
${k^1}_{j_1j_2}$ between the respective restrictions of
${e^1}_{j_2}$ and ${e^1}_{j_1}$ to ${U^0}_{j_1j_2}$.

 The map
 ${c^1}_{j_1j_2j_3}= {g^1}_{j_1j_2}{g^1}_{j_2j_3}{g^1}_{j_3j_1}$ is an
automorphism of   ${e^1}^{j_1}_{j_2j_3}$ the restriction of
$e^1_{j_1}$ to $U_{j_1j_2j_3}$ of  identified to
${U^0}_{j_1j_2j_3}\times F_1$. the chain $c_{j_1j_2j_3}$  can be
viewed as a map ${U^0}_{j_1j_2j_3}\rightarrow H_1$. Since $C_1$ is
a gerbe, the family ${c^1}_{j_1j_2j_3}$ is a $2-H_1$ cocycle.

The map ${c^1}_{j_1j_2j_3}$ acts on the category of open sets of
${U^0}_{j_1j_2j_3}\times F_1$.

The covering which defines the trivialization of
${e^1}_{j_1j_2j_3}$ can be considered as ${U^0}_{j_1j_2j_3}\times
{U^1}$. The map $c_{j_1j_2j_3}$ induces a functor
$(c_{j_1j_2j_3})^*$ between the categories
$C_2(U^0_{j_1j_2j_3}\times U^1)$ and $C_2(U^0_{j_1j_2j_3}\times
c_{j_1j_2j_3}(U^1))$. Since we have assumed $H_1$ central in
$G_1$, we can consider $c_{j_2j_3j_4}$, $c_{j_1j_3j_4}$,
$c_{j_1j_2j_4}$ as morphisms of $U^0_{j_1j_2j_3}\times F_1$, this
enable to compose the functor $(c_{j_1j_2j_3})*$
$(c_{j_2j_3j_4})^*$, $(c_{j_1j_3j_4})^*$, $(c_{j_1j_2j_4})^*$, and
to define the Cech boundary $d((c_{j_1j_2j_3})^*)$.

Since $c_{j_1j_2j_3}$ is a $2-$cocycle, the map
$d((c_{j_1j_2j_3{j^1}})^*)$ coincide with the restriction of an
element $c_{j_1j_2j_3j_4}$ of the set of differentiable functions
$U^0_{j_1j_2j_3j_4}\times U^1\longrightarrow H_2$ wich acts by
automorphisms on the elements of the category
$C_2(U_{j_1j_2j_3J_4}\times U^1)$.

\medskip

{\bf Proposition 1.}

{\it The family $c_{j_1j_2j_3j_4}$ defines an $3-$Cech cocycle.}

\medskip

{\bf Proof.}

We have $d(c_{j_1j_2j_3j_4})=\sum_{i=1}^{i=5} (-1)^ic_{j_1..\hat
j_i..j_5}=\sum_{i=1}^{i=5}(-1)^i\sum_{k=1}^{k=3}(-1)^k(c_{j_1..\hat
j_k..\hat j_i..})^*=0$.

\bigskip

Supposed now define the $i+1$ cocycle associated to the
$i-$lifting problem. It is a family of maps
$c_{j_1..j_{i+2}}:{U^0}_{{j}_1..{j}_{i+2}}\times
..\times{U^{i-1}}\times F_i$. The map $c_{j_1..j_{i+2}}$ defines a
functor $c(_{j_1..j_{i+2}})^*$ between the categories
$C_{i+1}(U_{j_1..j_{i+2}}\times U^1\times..\times U^i)$, and the
category $C_{i+1}(U_{j_1..j_{i+2}}\times U^1\times..\times U^i
c_{j_1..j_{i+2}}(U^i))$, where $U^i$ is an open set of $F_i$.
Since we have supposed that $H_i$ is central in $G_i$, we can
consider the Cech boundary $d(c_{j_1..j_{i+2}})^*$. Since the map
$c_{j_1..j_{i+2}}$ is an $i-$cocycle, we can consider
$d(c_{j_1..j_{i+3}})^*$ as a differentiable map
$c_{j_1..j_{i+3}}:U_{j_1..j_{i+3}}\times U^1..\times
U^i\rightarrow H_{i+1}$.

\bigskip

{\bf Proposition 2.} {\it The family of maps $c_{j_1..j_{i+3}}$
define an $j+2$ cocycle.}

{\medskip}

{\bf Proof.}

We have
$$
d(c_{j_1..j_{i+3}})=\sum_{k=1}^{k=i+4}(-1)^{k}c_{j_1.\hat j_k
j_{i+4}}=\sum_{k=1}^{k=i+4}(-1)^k\sum_{l=1}^{l=i+3}(-1)^l(c_{i_1..\hat
j_l..\hat j_k..j_{i+4}})^*=0
$$

\medskip

We have seen how to each $n-$lifting problem one can associate a
$n+1-$cocycle $H_n$ cocycle, is it true that fact conversely to
each $n+1-H_n$, one can define a $n+1-$ lifting problem which
$n+1-$ classifying cocycle is the given cocycle? if not which
conditions must satisfy and $n+1-$cocycle to be the classifying
cocycle of a $n+1-$lifting problem? We will first study the first
lifting problem. Recall that if we have the following exact
sequence $1\rightarrow S_1\rightarrow S_2\rightarrow
S_3\rightarrow 1$ of sheaves defined on $B$ such that $S_1$ is a
commutative sheaf, there exists the following exact sequence in
cohomology.
$$
0\rightarrow H^0(B,S_1)\rightarrow H^0(B,S_2)\rightarrow
H^0(B,S_3)
$$
$$
\rightarrow H^1(B,S_1)\rightarrow H^1(B,S_2)\rightarrow
H^1(B,S_3)\rightarrow H^2(B,S_3)
$$
denote by $\delta_1$ the map $H^1(B,S_3)\rightarrow H^2(B,S_3)$,
we have

\medskip

{\bf Proposition 3.}

{\it A $2-$$H_1$cocycle $c_{H_1}$ defined on $B$ is the classying
cocycle of a first lifting problem, if and only if there exists an
element $c_{K_1}$ of $H^1(B,K_1)$ such that
$c_{H_1}=\delta_1(c_{K_1})$.}

\medskip

{\bf Proof.}

Its enough to remark that each $c_{K_1}$ cocycle $c$ define a
principal bundle over $B$, such that $\delta_1(c)$ is the
classifying cocycle of the first lifting problem associated to
this bundle.

\medskip

{\bf Remark.}

Saying that $\delta_1(c_{K_1})=c_{H_1}$ means that $c_{H_1}$  is
deduced from $c_{K_1}$ in the same way that $\delta_1$ is defined.

\medskip

Generally the map $\delta_1$ is not injective. To see this
consider a $G-$principal bundle $p_{G}$ over $B$, $H$ a central
subgroup of $G$. Denote by $K$ the quotient of $G$ by $H$ and by
$c$ the cocycle which defines the quotient $p_{K}$ of $p_{G}$ by
$H$, then $\delta_1(c)$ is trivial.

\medskip

Now suppose that the sheaves $S_1$, $S_2$ and $S_3$ are
commutative, we can write the following long exact sequence

$$
0\rightarrow H^0(B,S_1)\rightarrow H^0(B,S_2)\rightarrow
H^0(B,S_3)
$$
$$
\rightarrow H^1(B,S_1)\rightarrow H^1(B,S_2)\rightarrow
H^1(B,S_3)\rightarrow H^2(B,S_1)
$$
$$
..\rightarrow H^n(B,S_1)\rightarrow H^n(B,S_2)\rightarrow H^n
(B,S_3)\rightarrow H^n(B,S_3)
$$
denote by $\delta_n$ the map $H^n(B,S_3)\rightarrow
H^{n+1}(B,S_3)$, we have

\medskip

{\bf Proposition 4.}

{\it The $n+1-$ cocycle $c_{H_n}$ of the $n-$lifting problem is
the image of $c_{H_{n-1}}$ by $\delta_n$}

\medskip

{\bf Proof.}

Straightforward computations.

\medskip

{\bf Definition 1.}

We will say that the $n+1-$tower is trivial, if and only if its
associated cocycle $c_{H_n}$ vanishes.

\bigskip

{\bf The principal tower of gerbes associated.}

\bigskip

We will define the notion of connective structure and curvature
related to the tower of gerbes associated to our problem. First we
define the principal tower associated to our problem.

\medskip

Recall that we can associate to $C_1$ a principal gerbe $C_{1P}$
defined as follows:

Let $U$ be an open set of $B$, we consider the category
$C_{1P}(U)$ of $G_1$ principal bundles over $U$ such that the
quotient of each of its elements by $H_1$ is the restriction of
$p_{K_1}$ to $U$,

Let $e_1(U)$ be an element of $C_1(U)$, we consider an open set
$U^1$ of $e_1(U)$, and an element $e_2$ of $C_2(U^1)$, it is an
$F_2$ bundle over $U^1$, we denote by $p_{U^1}:e_2\rightarrow
U^1$, the canonical projection, we denote by
$p_{1U}:U^1\rightarrow U$ the restriction of the canonical
projection of $e_1(U)$ to $U^1$.

We can also defined over $e_2$ the principal bundle $e_{2P}$
associated to the first lifting problem defined by $F_2$, $H_2$
and $G_2$ which typical fiber is $G_2$, the elements of the
category $C_{2P}(U^1)$, are fiber product over $p_{1U}(U^1)$ of
$e_{2P}$ and $e_{1P}$, where $e_{1P}$ is a principal bundle over
$p_{1U}(U^1)$ associated to the first lifting problem defined by
the data $F_1$, $H_1$ and $G_1$. Suppose that $U^1=V\times V^1$
where $V$ and $V^1$ are open sets of respectively $B$ and $F_1$
such that the respective restriction of $p_{K_1}$ and $p_{K_2}$ to
$V$ and $V^1$ are trivial. Then the  elements of $C_{2P}(U^1)$ are
isomorphic to $(V\times G_1)\times (V_1\times G_2)$. The set of
morphisms of elements of $C_{2P}(U^1)$ is given by the natural
action of $H_2$.

Suppose defined the category $C_{iP}(U^{i-1})$, where $U^{i-1}$ is
an open set of an element $C_{i-1}(U^{i-2})$, let $U^{i}$ be an
open set of an object of $C_i(U^{i-1})$, we can define on $U^i$
the principal bundle $e_{i+1P}$ which solve the first lifting
problem associated to $F_{i+1}$, $G_{i+1}$ and $H_{i+1}$. It is a
principal bundle $e_{iP}$ over $U^i$ which typical fiber is
$G_{i+1}$, denote by $p_{iU}$ the canonical map
$p_{iU}:U^i\rightarrow U^{i-1}$, the elements of $C_{i+1P}(U^i)$,
are fiber product of $e_{iP}$ and $e_{i-1P}$ over
$p_{iU}(U^{i-1})$ where $e_{i-1P}$ is an element of
$C_{iP}(p_{iU}(U^{i}))$. The automorphisms group of elements of
$C_{i+1P}(U^i)$ are induced by the action of $H_{i+1}$.

\medskip

Locally, an object of $C_{iP}(U^{i-1})$ is of the form $U\times
G_1\times(U^1\times G_2)\times..\times (U^{i-1}\times G_i)$.

\medskip

{\bf Definition 1.}

The $n-$tangent space of the tower of gerbes $C_P$, will be the
disjoint unions of the tangent space of the gerbes $C_{iP}$.

\bigskip

{\bf The differential geometry of the principal tower of gerbes.}

\medskip

Suppose defined on each bundle $p_{K_i}$ a connection $w_i$. We
have already define the notion of connective structure for the
first lifting problem. Recall that for each element $e_1$ of
$C_{1P}(U)$, it is the family of connections defined on $e_1$
which project to the restriction of $w_1$ to the restriction of
$p_{K_1}$ to $U$.

Suppose defined the notion of connective structure for the
$i-$lifting problem, for each object $e_{iP}$ of
$C_{iP}(U^{i-1})$, (we can suppose $e_{iP}$ to be isomorphic to
$U\times G_1\times(U^1\times G_2)\times..\times (U^{i-1}\times
G_i)$)
  it is a family
of connection defined on $e_{iP}$, for an element $e_{i+1P}$ of
$C_{i+1P}(U^i)$, we will defined the the set $Co(e_{i+1P})$, as
the set of $G_1\times..\times G_{i+1}$connections  which project
to elements of $Co(e_i(p_{iU}(U^i))$, where $e_i(p_{U_i}(U^i))$ is
the canonical projection of $e_{i+1P}$ over an element of
$C_{iP}((p_{iU})(U_i))$, consider the projection of $U^i$ to
$F_i$, and the restriction $p_{UK_i}$ of $p_{K_i}$ to the image of
the last projection, we suppose also that  elements of
$Co(e_{i+1P})$ project to the restriction of $w_{i+1}$.

If this element is isomorphic to $U\times G_1\times..\times
U^{i-1}\times G_i$, then an element of $Co(e_{iP})$ is a
connection which project to a connection of $Co(e_{i-1P})$, where
$e_{i-1P}=U\times G_1\times..\times U^{i-2}\times G_{i-2}$ and
which project to a connection of $U^{i-1}\times G_i$, which
project to the restriction of $w_i$ to $U^{i-1}\times K_i$.

\bigskip

{\bf The sequence of curvatures associated.}

\medskip

For the gerbe associated to the first lifting problem, we have
already defined the curvature, it is a $Lie(1)$ $3-$form. Let us
recall it's definition.

Consider a trivialization $(U_{i_1},k_{i_1i_2})$ of the bundle
$p_{K_1}$, over $U_{i_1}\times G_1$ we choose a connection
$w_{i_1}$ which project to the restriction of $w_1$ to
$U_{i_1}\times G_1$, we consider the curvature $K(w_{i_1})$ of
this connection, then $dK(w_{i_1})$ defines the requested
$3-$form.

Consider now an element $e_i$ of $C_{iP}(U^{i-1})$, isomorphic to
$U_{j_1..j_{i+3}}\times G_1\times(U^1\times G_2)\times..\times
(U^{i-1}\times G_i)$, and $\hat w_i$ an element of $Co(e_i)$, if
$\hat w'_i$ is another element of $Co(e_i)$, then $\hat w_i-\hat
w'_i$ is an $H_1\oplus..\oplus H_n$ form $a$, we have
$K(w_i)=K(w'_i)+da$, this implies that $d(K(w_i))$ is independent
of the chosen element of $Co(e_i)$ over $w_i$

\medskip

{\bf Definition 2.}

 The form curvature $dK(\hat w_i)$, defines
a $3-$form $\Omega_{i}$ on $3TC_{iP}(U^{i-1})$.

\bigskip

{\bf III. Holonomy and parallel displacement.}

\bigskip

We will precise the differential geometry of the $n-$lifting
problem   by defining notions of holonomy and covariant
derivative.

Suppose given a connection $\omega$ on $p_{K_1}$, for each open
set $U$ of $B$ and each couple of objects $(e_1,e_2)$ of objects
of $C_1(U)$, one can consider connections $w_1$ and $w_2$ defined
respectively on $e_1$ and $e_2$ which push forward to the
restriction of $\omega$ to the restriction of $p_{K_1}$ to $U$,
unfortunately each isomorphism between $e_1$ and $e_2$ does not
transforms $\omega_1$ to $\omega_2$. To avoid this problem, we are
going to enlarge the distribution which defines the connection
$\omega_1$ and $\omega_2$, to a transitive uniform distribution
studied by P. Molino in his thesis, this distribution will behave
naturally in respect to the morphisms between objects,
generalizing the work of the previous author, will define first,
the notion of holonomy and covariant derivative for the principal
gerbe $C_1$ and generalizing those notions for the principal tower
associated to the last $n-$lifting problem. First let recall some
results from the work of P. Molino.

\medskip

Let $p:E\rightarrow B$, be a $K-$principal locally trivial bundle
over $B$ which typical fiber is $K$, a transitive distribution
${\cal D}$ (TD) on $E$ will be a right invariant distribution
defined on $E$ such that $d_xp({\cal D}_x)=T_{p(x)}B$.

We will denote by $u$ be the map defined by $u(X)=A$, where $X$ a
vertical vector tangent to $x$ and $A$ the element of ${\cal K}$
such that $A^*(x)=X$, where $A^*$ is  the fundamental vector field
defined by $A$.    We will also denote $V_x$ the subset of
vertical vector fields at $x$. The transitive distribution is said
to be freely uniform, if $u({\cal D}\cap V_x)$ does not depend of
$x$, this image is thus a normal subgroup $L$ of $K$.

In the sequel, we will only consider free transitive uniform
distributions.

\medskip

{\bf Definition 1.}

An horizontal curve of ${\cal D}$ is a curve $c:I\rightarrow E$
such that, $c'(t)={d\over dt}c(t)$ is an element of ${\cal D
}_{c(t)}.$

\medskip

Let, $x$ be an element of $E$, the holonomy group of ${\cal D}$ at
$x$, is the set of elements $k$ of $K$ such that $x$ and $xk^{-1}$
can be joined by an horizontal curve.

\medskip

Let denote by $M$ the quotient of $K$ by $L$, $U:K\rightarrow M$
the projection map, and $d:{\cal K}\rightarrow {\cal M}$ the
induced Lie algebra morphism, one can define on $E$, a $1-$ $M$
form as follows, each element $v$ of $T_xE$ can be written
$v=v_1+v_2$, where $v_1$ is vertical, and $v_2$ is an element of
${\cal D}_x$, this decomposition is not unique unless that ${\cal
D}$ is a connection, but the projection $\omega(v)$ of $u(v_1)$ to
${\cal M}$ defined on $E$ a $1-$form which kernel is ${\cal D}$.

\medskip

Let $V$ be a vector space endowed with a $K-$action defined by the
representation $S:K\rightarrow Gl(V)$, we will denote by $s:{\cal
K}\rightarrow gl(V)$ the induced Lie algebras representation,  by
$V'$ the quotient of $V$ by $s({\cal L})$, and by $s':V\rightarrow
V'$ the projection map where ${\cal L}$ is the Lie algebra of $L$.
We consider a $p-$tensorial form $\alpha: TE\rightarrow V$. We
denote $\dot\alpha$, the $V'$ $p-$tensorial form defined by
$\dot\alpha(v_1,..,v_p)=s'(\alpha(v_1,...,v_p))$ where
$v_1,..,v_p$ are elements of $T_xE$, write $v_i=w^i_1+w^i_2$,
where $w^i_1$ is an vertical vector and $w^i_2$ is an element of
${\cal D}_x$ we will define the covariant derivative,
$\nabla_{\cal D}$ of $\dot\alpha$ by $\nabla_{\cal
D}(v_1,...,v_p)=s'(d\dot\alpha(w^1_2,...,w^p_2))$.

We have the formula $\nabla_{\cal
D}(\dot\alpha)=d\dot\alpha+\omega\dot\alpha$. We will denote by
$\nabla_{D}\dot\omega$, the curvature of ${\cal D}$.

\bigskip

We will now apply this notion of TD to the tower of gerbe. First
we consider the case of the gerbe $C_1$. Let endow $p_{K_1}$ with
a connection $w_1$ and consider an object $e_1$ of $C_1(U)$, where
$U$ is an open set of $B$. We can define on $e_1$ the transitive
distribution ${{\cal D}_1}_{e_1}$, the pull-back of the
distribution which define the connection $w_1$ by the canonical
map $e_1\rightarrow {p_{K_1}}_{\mid U}$, if $e_1'$ is an other
object of $C(U)$, then each map $f:e_1\rightarrow e_1'$ transforms
${{\cal D}_1}_{e_1}$ to ${{\cal D}_1}_{e_1'}$. The TD ${\cal
D}_{e_1}$ is uniform, its uniform group is $H_1$.

\medskip

{\bf Definition 2.}

For each element $z$ of $e_1$, we will denote by
$Hol(z,e_1,\omega,w_1)$, the holonomy group of the transitive
distribution ${\cal D}_{e_1}$ in $z$. Since we have supposed $H_1$
to be central in $G_1$, this group depends only of the projection
of $z$ to $p_{K_1}$.

\bigskip

{\bf The holonomy of the tower.}

\bigskip

We will define recursively the notion of holonomy of the tower of
gerbe, by defining on each element of $C_{i+1P}(U^i)$ a free
transitive distribution. We suppose defined on each bundle
$p_{K_i}$ a connection $w_i$.

\medskip

Let $U^1$ be an open set of $e_1$, recall and object $e^2_{U^1}$
of $C_{2P}(U^1)$, are fibers product over $p_{1U}(U^1)$, of
$e_{2P}$ and $e_{1P}$, where $e_{2P}$ is a principal bundle
associated to an object of $C_2(U^1)$ associated to the first
lifting problem defined by $H_2$, $F_2$ and $G_2$, and $e_{1P}$ is
the principal bundle over $p_{1U}(U^1)$ of the first lifting
problem associated to the data $H_1$, $F_1$ and $G_1$.

We have two canonical projections $f_1: e^2_{U^1}\rightarrow
e_{1P}$, and $f_2: e^2_{U^1}\rightarrow e_{2P}$.

The connection $w_2$ on $p_{K_2}$ can be pull back to a connection
$w'_2$ on the quotient of $e_{2P}$ by $H_2$.
 Denote by
$X_{e_{2P}}$ the pull back to $e_{2P}$ of the distribution which
defines the connection $w'_2$ to $e_{2P}$, we will define
$$
{{\cal D}_2}_{e^2_{U^1}}=f_1^* {{\cal
D}_1}_{e_{1P}}+f_2^*X_{e_{2P}}.
$$

\bigskip

{\bf Definition 3.}

 Let $z$ be an element of $e^2_{U^1}$. We denote
by $Hol(e^2_{U^1},{\cal D}_2,z, w_1,w_2)$ the holonomy of the
transitive distribution ${{\cal D}_2}_{e^2_{U^1}}$ in $z$.

\medskip

Suppose defined the distribution ${{{\cal D}_i}_{e_{iP}}}$ and the
holonomy $Hol(e_{iP},{\cal D}_i,z,w_1,..,w_i)$ where is an element
of the object $e_{iP}$, of $C_{iP}(U^{i-1})$.

 Recall that an object
$e^{i+1}$ of $C_{i+1P}(U^i)$ is a fiber product over $p_{iU}(U^i)$
of $e_{iP}$ and $e_{i+1P}$, where $e_{iP}$ is an element of
$C_{iP}(U^{i-1})$, and $e_{i+1P}$ is a bundle over $U^i$ which
solve the first lifting problem associated to $F_{i+1}$, $H_{i+1}$
and $G_{i+1}$, we have the canonical projection
$f_1:e^{i+1}\rightarrow e_{iP}$, and $f_2:e^{i+1}\rightarrow
e_{i+1P}$. The connection $w_{i+1}$ on $p_{K_{i+1}}$ induces on
the quotient of $e_{i+1P}$ by $H_{i+1}$ a connection
$\omega'_{i+1}$. Denote by $X_{i+1}$ the pull back  of the
distribution which defines the connection $\omega'_{i+1}$ to
$e_{i+1P}$, we can define

$$
{{\cal D}_{i+1}}_{e^{i+1}}=f_1^*({{\cal
D}_i}_{e_{iP}})+f_2^*(X_{i+1}).
$$

\bigskip

{\bf Definition 4.}

For each element $z$ of $e^{i+1}$, we denote by $Hol(e^{i+1},{\cal
D}_{i+1},z,w_1,..,w_{i+1})$, the holonomy of the  transitive
distribution ${{\cal D}_{i+1}}_{e^{i+1}}$ in $z$.

\bigskip

{\bf Parallel Displacement.}

\medskip

 We will define the notion of parallel displacement for the tower
 of gerbes.

 \medskip

 {\bf Definition 5.}

  Given a connection $w_1$ on the principal bundle
 $p_{K_1}$, and $z$ $c:I=[0,1]\rightarrow p_{K_1}$ a differentiable curve,
A field $X$ is parallel along $c$, if and only if
${\nabla_{w_1}}_{{d\over dt}c(t)}X(c(t))=0$, where $\nabla_{w_1}$
is the covariant derivative associated to the connection $w_1$.

 The parallel
 displacement along $c$, is the map $s_{c(t)}:T_{c(0)p_{K_1}}
 \rightarrow T_{c(1)}p_{K_1}$
such that $s_{c(t)}(v_0)$ is the vector $v_1$ of $T_{c(1)}p_{K_1}$
such that $v_0=X({c(0)}), v_1=X(c(1))$ where $X$ is a parallel
vector field defined on $c$, (we suppose $c$ endowed the pull-back
of the connection $w_1$).

\bigskip

{\bf Definition 6.}

Consider an object $e_{k}$ of the category $C_{kP}(U^k)$ of the
principal tower, and a connection $\omega$ of $Co(e_k)$. For each
element $z$ of $e_k$, for each curve $c:I=[0,1]\rightarrow e_k$,
one can define the parallel displacement $s_{c,\omega}$ along $c$,
the parallel displacement along $c$ is the set of maps
$s_{c,\omega}, \omega\in Co(e_k)$.

\bigskip

\centerline{\bf Bibliography.}

\bigskip

[1]  Brylinski, J.L Loops spaces, Characteristic Classes and
Geometric Quantization, Progr. Math. 107, Birkhauser, 1993.

[2] Brylinski, J.L, Mc Laughlin D.A, The geometry of degree four
characteristic classes and of line bundles on loop spaces I. Duke
Math. Journal. 75 (1994) 603-637.

[3] Deligne, P. Theorie de Hod ge III, Inst. Hautes Etudes Sci.
Publ. Math. 44 (1974), $5-77$.

[4] Giraud, J. Cohomologie non ab\'elienne. Grundl. Math. Wiss.
179, Springer-Verlag.

[5] Goettsche, L. Hilbert Schemes. Lecture note in Mathematics

[6] Meinrenken, E. The basic gerbe of a compact simple Lie group.
Preprint

[7] Molino, P. Champs d'\'el\'ements sur un espace fibr\'e
principal diff\'erentiable. Ann. Inst. Fourier 14 (1964) 163-220.

[8] Murray, M.K., Stevenson, Bundle gerbes: Stable isomorphisms
and local theory. J. London. Math. Soc (2) 62 (2000) 925-937

[9] Perea Medina, A. Flat left-invariant connections adapted to
the automorphism structure of a Lie group, J. Differential
Geometry, 16 (1981) 445-474.

[10] Tsemo, A. Non abelian cohomology an application Preprint.

\end{document}